\documentclass{amsart}
\usepackage{latexsym}
\usepackage{amsfonts}
\usepackage{amssymb}
\usepackage{graphicx}
\usepackage{graphicx}

\newtheorem{theorem}{Theorem}[section]

\newtheorem{lemma}[theorem]{Lemma}
\newtheorem{proposition}[theorem]{Proposition}
\theoremstyle{definition}

\theoremstyle{remark}
\newtheorem{remark}{Remark}[section]

\newcommand{\g}{\mathfrak{g}}
\newcommand{\ak}{\mathfrak{k}}
\newcommand{\ap}{\mathfrak{p}}
\newcommand{\al}{\mathfrak{l}}
\newcommand{\am}{\mathfrak{m}}
\newcommand{\ee}{\textbf{e}}

\begin{document}

\title[]
{On the realization of Riemannian symmetric spaces in Lie groups
II$^1$}

\author[]{Jinpeng An}
\address{School of mathematical science, Peking University,
 Beijing, 100871, P. R. China }
\email{anjinpeng@gmail.com}

\author[]{Zhengdong Wang}
\address{School of mathematical science, Peking University,
 Beijing, 100871, P. R. China}
\email{zdwang@pku.edu.cn}

\maketitle

{\bf Abstract.} In this paper we generalize a result in \cite{AW},
showing that an arbitrary Riemannian symmetric space can be
realized as a closed submanifold of a covering group of the Lie
group defining the symmetric space. Some properties of the
subgroups of fixed points of
involutions are also proved.\\


\footnotetext[1]{This work is supported by the 973 Project
Foundation of China ($\sharp$TG1999075102).}


\vskip 1.0cm
\section{Introduction}

\vskip 0.3cm Suppose $G$ is a connected Lie group with an
involution $\sigma$. Then the Lie algebra $\g$ of $G$ has a
canonical decomposition $\g=\ak\oplus\ap$, where $\ak$ and $\ap$
are the eigenspaces of $d\sigma$ in $\g$ with eigenvalues $1$ and
$-1$, respectively. Let $G^\sigma=\{g\in G|\sigma(g)=g\}$, and
suppose $K$ is an open subgroup of $G^\sigma$. Suppose moreover
that $Ad_G(K)|_\ap$ is compact, that is, $G/K$ has a structure of
Riemannian symmetric space. In the particular case that
$K=G^\sigma$, it was proved in \cite{AW} that $P=\exp(\ap)$ is a
closed submanifold of $G$, and there is a natural isomorphism
$G/G^\sigma\cong P$. This gave a realization of the symmetric
space $G/G^\sigma$ in $G$. In this paper we generalize this result
to the case of arbitrary symmetric space $G/K$, that is, the case
that $K$ is an arbitrary open subgroup of $G^\sigma$ such that
$Ad_G(K)|_\ap$ is compact.

\vskip 0.3cm In Section 2 we will make some preparation for this
generalized realization. Some properties of the subgroup
$G^\sigma$ will be examined. We will prove the following results.

\begin{itemize}
\item{There exists a covering group $G'$ of $G$ with covering
homomorphism $\pi$ such that there is an involution $\sigma'$ on
$G'$ with $\pi\circ\sigma'=\sigma\circ\pi$, and such that
$\pi^{-1}(K)=(G')^{\sigma'}$.}

\item{The quotient group $G^\sigma/G^\sigma_0$ is isomorphic to
$(\mathbb{Z}_2)^r$ for some non-negative integer $r$, where
$G^\sigma_0$ is the identity component of $G^\sigma$.}

\item{$G^\sigma$ is connected if $\pi_1(G)$ is finite with odd
order.}
\end{itemize}

\vskip 0.3cm In Section 3 we will give the precise statement of
the realization of arbitrary symmetric spaces in Lie groups.
Briefly speaking, a symmetric space $G/K$ is diffeomorphic to a
closed submanifold $P'$ of a covering group $G'$ of $G$, where
$G'$ is chosen such that $\pi^{-1}(K)=(G')^{\sigma'}$. The idea of
the proof is that $G/K\cong G'/\pi^{-1}(K)=G'/(G')^{\sigma'}\cong
\exp_{G'}(\ap)=P'$.

\vskip 0.3cm The authors would like to thank the referee for some
suggestion on improving the English usage of the paper. They also
thank the referee of \cite{AW}, whose suggestion of revision
motivated the authors considering the problems of the paper. The
first author thank Professor Jiu-Kang Yu for valuable
conversation.


\vskip 1.0cm
\section{The subgroups of fixed points of involutions}

\vskip 0.3cm For a Lie group $H$ with an automorphism $\theta$, we
always denote $H^\theta=\{h\in H|\theta(h)=h\}$, and denote the
identity component of $H^\theta$ by $H^\theta_0$. In this section
we prove the following two theorems. The realization of arbitrary
symmetric spaces in Lie groups will be based on Theorem
\ref{T:covering}.

\begin{theorem}\label{T:covering}
Let $G$ be a connected Lie group with an involution $\sigma$, K an
open subgroup of $G^\sigma$. Then there exists a covering group
$G'$ of $G$ with covering homomorphism $\pi$ such that there is an
involution $\sigma'$ on $G'$ with
$\pi\circ\sigma'=\sigma\circ\pi$, and such that
$\pi^{-1}(K)=(G')^{\sigma'}$.
\end{theorem}

\begin{theorem}\label{T:quotient}
Let $G$ be a connected Lie group with an involution $\sigma$. Then
the quotient group $G^\sigma/G^\sigma_0$ is isomorphic to
$(\mathbb{Z}_2)^r$ for some non-negative integer $r$, and $r=0$ if
$\pi_1(G)$ is finite with odd order.
\end{theorem}

First we introduce some notations. Let $G$ be a connected Lie
group with an involution $\sigma$. For $g,h\in G$, denote the set
of all continuous paths $\gamma:[0,1]\rightarrow G$ with
$\gamma(0)=g$ and $\gamma(1)=h$ by $\Omega(G,g,h)$, and denote
$\pi_1(G,g,h)=\{[\gamma]|\gamma\in\Omega(G,g,h)\}$, where
$[\gamma]$ is the homotopy class relative to endpoints determined
by $\gamma$. Note that $\pi_1(G,g,g)$ is the fundamental group
$\pi_1(G,g)$ of $G$ with basepoint $g$. For $g\in G^\sigma$ and
$\gamma_g\in\Omega(G,e,g)$, let $\gamma_g^\sigma\in\Omega(G,e,e)$
be defined as
$$
\gamma_g^\sigma(t)=
\begin{cases}
\gamma_g(2t) & t\in[0,\frac{1}{2}];\\
\sigma(\gamma_g(2-2t)) & t\in[\frac{1}{2},1].
\end{cases}
$$
The set $\pi_g=\{[\gamma_g^\sigma]|\gamma_g\in\Omega(G,e,g)\}$ is
a subset of $\pi_1(G,e)$. For $g_1,g_2\in G^\sigma$, if they
belong to the same coset space of $G^\sigma_0$, it is obvious that
$\pi_{g_1}=\pi_{g_2}$. So we can define $\pi_{[g]}=\pi_g$ for
$g\in G^\sigma$, where $[g]$ denotes the coset space
$gG^\sigma_0$. We denote the identity element of $\pi_1(G,e)$ by
$\ee$.

\begin{lemma}\label{L:inverseclosed}
For $g\in G^\sigma$, $\pi_{[g]}^{-1}=\pi_{[g]}$, that is,
$x\in\pi_{[g]}\Leftrightarrow x^{-1}\in\pi_{[g]}$.
\end{lemma}

\begin{proof}
It is obvious from the equation
$[\gamma_g^\sigma]^{-1}=[(\sigma\circ\gamma_g)^\sigma]$.
\end{proof}

\begin{lemma}\label{L:order2}
For $g\in G^\sigma$, $\pi_{[g]}=\pi_{[g^{-1}]}$.
\end{lemma}

\begin{proof}
For $[\gamma_g^\sigma]\in\pi_{[g]}$,
$[(\gamma_g^{-1})^\sigma]\in\pi_{[g^{-1}]}$, where
$\gamma_g^{-1}(t)=\gamma_g(t)^{-1}$. But by Lemma 16.7 in
\cite{St}, $[\gamma_g^\sigma]\cdot[(\gamma_g^{-1})^\sigma]
=[\gamma_g^\sigma]\cdot[(\gamma_g^\sigma)^{-1}]
=[\gamma_g^\sigma(\gamma_g^\sigma)^{-1}]=\ee$, where
$(\gamma_g^\sigma(\gamma_g^\sigma)^{-1})(t)=\gamma_g^\sigma(t)(\gamma_g^\sigma)^{-1}(t)$.
So by Lemma \ref{L:inverseclosed},
$[\gamma_g^\sigma]=[(\gamma_g^{-1})^\sigma]^{-1}\in\pi_{[g^{-1}]}^{-1}=\pi_{[g^{-1}]}$.
This proves $\pi_{[g]}\subset\pi_{[g^{-1}]}$. By symmetry, we have
the equality.
\end{proof}

\begin{lemma}\label{L:homomorphism}
For $g_1,g_2\in G^\sigma$,
$[\gamma_{g_1}^\sigma]\cdot\pi_{[g_2]}=\pi_{[g_1g_2]}$ for each
$[\gamma_{g_1}^\sigma]\in\pi_{[g_1]}$. Hence
$\pi_{[g_1]}\cdot\pi_{[g_2]}=\pi_{[g_1g_2]}$.
\end{lemma}

\begin{proof}
Let $[\gamma_{g_2}^\sigma]\in\pi_{[g_2]}$. By Lemma 16.7 in
\cite{St}, $[\gamma_{g_1}^\sigma]\cdot[\gamma_{g_2}^\sigma]
=[\gamma_{g_1}^\sigma\gamma_{g_2}^\sigma]=[(\gamma_{g_1}\gamma_{g_2})^\sigma]
\in\pi_{[g_1g_2]}$. So we have
$[\gamma_{g_1}^\sigma]\cdot\pi_{[g_2]}\subset\pi_{[g_1g_2]}$. To
prove the equality, denote $x=[\gamma_{g_1}^\sigma]$, and consider
the map $l_x:\pi_{[g_2]}\rightarrow\pi_{[g_1g_2]}$ defined by
$l_x(y)=xy$. Since
$x^{-1}\in\pi_{[g_1]}^{-1}=\pi_{[g_1]}=\pi_{[g_1^{-1}]}$,
$x^{-1}\cdot\pi_{[g_1g_2]}\subset\pi_{[g_2]}$. So we can also
define the map $l_{x^{-1}}:\pi_{[g_1g_2]}\rightarrow\pi_{[g_2]}$
by $l_{x^{-1}}(y)=x^{-1}y$. Since $l_x\circ l_{x^{-1}}=id$, $l_x$
is surjective. This means $x\cdot\pi_{[g_2]}=\pi_{[g_1g_2]}$,
hence $\pi_{[g_1]}\cdot\pi_{[g_2]}=\pi_{[g_1g_2]}$.
\end{proof}

\begin{lemma}\label{L:piK}
Let $K$ be an open subgroup of $G^\sigma$. Then
$\pi_K:=\bigcup_{k\in K}\pi_{[k]}$ is a subgroup of $\pi_1(G,e)$.
\end{lemma}

\begin{proof}
By Lemma \ref{L:inverseclosed}, each $\pi_{[k]}$ is closed under
the inverse operation, hence so is $\pi_K$. By Lemma
\ref{L:homomorphism}, $\pi_K$ is also closed under multiplication.
\end{proof}

By Lemma \ref{L:piK}, $\pi_{[e]}=\pi_{(G^\sigma_0)}$ and
$\pi_{(G^\sigma)}$ are subgroups of $\pi_1(G,e)$. By Lemma
\ref{L:homomorphism}, for each $g\in G^\sigma$ and each
$x\in\pi_{[g]}$, $\pi_{[g]}=x\cdot\pi_{[e]}$. So $\pi_{[g]}$ is a
coset space of $\pi_{[e]}$ in $\pi_{(G^\sigma)}$, that is, an
element in the quotient group
$\Pi^\sigma=\pi_{(G^\sigma)}/\pi_{[e]}$ (note that the fundamental
group of a Lie group is always abelian). Define the map
$f:G^\sigma/G^\sigma_0\rightarrow\Pi^\sigma$ by
$f([g])=\pi_{[g]}$. By Lemma \ref{L:homomorphism}, $f$ is a
homomorphism.

\begin{lemma}\label{L:finite}
Let $G$ be a connected Lie group with an involution $\sigma$. Then
$G^\sigma$ has finite many connected components. If moreover $G$
is simply connected, then $G^\sigma$ is connected.
\end{lemma}

\begin{proof}
The subgroup $\Sigma=\{id, \sigma\}$ of the automorphism group
$Aut(G)$ of $G$ has a natural action on $G$, which we also denote
by $\sigma$. We form the semidirect product
$\mathbf{G}=G\times_\sigma\Sigma$. Denote
$x=(e,\sigma)\in\mathbf{G}$. Since the identity component
$\mathbf{G}_0$ of $\mathbf{G}$ is naturally isomorphic to $G$ by
$i:(g,id)\mapsto g$, and $\sigma(i(y))=i(xyx^{-1}), \forall
y\in\mathbf{G}_0$, to prove the lemma, it is sufficient to show
that $\mathbf{G}_0^x=\{y\in\mathbf{G}_0|xyx^{-1}=y\}$ has finite
many connected components, and is connected when $\mathbf{G}_0$ is
simply connected.

Since $x=(e,\sigma)$ is an element of order two in $\mathbf{G}$,
it lies in some maximal compact subgroup $\mathbf{L}$ of
$\mathbf{G}$. By Theorem 3.1 of Chapter XV in \cite{Ho}, there
exist some linear subspaces $\am_1,\cdots,\am_k$ of the Lie
algebra $\g$ of $\mathbf{G}$ (which is isomorphic to the Lie
algebra of $G$) such
that\\
(1) $\g=\al\oplus\am_1\oplus\cdots\oplus\am_k$, where $\al$ is the
Lie algebra of $\mathbf{L}$;\\
(2) $Ad(l)(\am_i)=\am_i, \forall l\in \mathbf{L},
i\in\{1,\cdots,k\}$
(in particular, $Ad(x)(\am_i)=\am_i$);\\
(3) the map
$\psi:\mathbf{L}_0\times\am_1\times\cdots\times\am_k\rightarrow
\mathbf{G}_0$ defined by $\psi(l,X_1,\cdots,X_k)=le^{X_1}\cdots
e^{X_k}$ is a diffeomorphism, where $\mathbf{L}_0$ is the identity
component of $\mathbf{L}$.\\
Denote $\mathbf{L}_0^x=\{l\in\mathbf{L}_0|xlx^{-1}=l\}$,
$\am_i^x=\{X\in\am_i|Ad(x)(X)=X\}, i=1,\cdots,k$. We claim that
$\mathbf{G}_0^x=\psi(\mathbf{L}_0^x\times\am_1^x\times\cdots\times\am_k^x)$.
In fact, for $y\in\mathbf{G}_0^x$, write $y=le^{X_1}\cdots
e^{X_k}$, where $l\in\mathbf{L}_0, X_i\in\am_i$. Then
$y=xyx^{-1}=(xlx^{-1})e^{Ad(x)(X_1)}\cdots e^{Ad(x)(X_k)}$. By
(3), we have $xlx^{-1}=l$, $Ad(x)(X_i)=X_i, i=1,\cdots,k$. That
is, $l\in\mathbf{L}_0^x$, $X_i\in\am_i^x$. Hence we have
$\mathbf{G}_0^x\subset\psi(\mathbf{L}_0^x\times\am_1^x\times\cdots\times\am_k^x)$.
The other inclusion is obvious.

Since $\mathbf{L}_0$ is a connected compact Lie group,
$\mathbf{L}_0^x$ has finitely many connected components, hence so
does
$\mathbf{G}_0^x=\psi(\mathbf{L}_0^x\times\am_1^x\times\cdots\times\am_k^x)$.
If $\mathbf{G}_0$ is simply connected, by (3), $\mathbf{L}_0$ is
also simply connected. By Theorem 8.2 of Chapter VII in \cite{He},
$\mathbf{L}_0^x$ is connected, hence so is $\mathbf{G}_0^x$. This
proves the lemma.
\end{proof}

\begin{lemma}\label{L:injective}
For $g\in G^\sigma$, if $\ee\in\pi_{[g]}$, then $[g]=G^\sigma_0$.
\end{lemma}

\begin{proof}
We endow $\widetilde{G}=\bigcup_{g\in G}\pi_1(G,e,g)$ with the
canonical smooth manifold structure and the canonical group
structure such that $\widetilde{G}$ is the universal covering
group of $G$ with covering homomorphism $\pi([\gamma])=\gamma(1)$.
The induced involution of $\sigma$ on $\widetilde{G}$ is
$\widetilde{\sigma}([\gamma])=[\sigma\circ\gamma]$. Let $g\in
G^\sigma$. If $\ee\in\pi_{[g]}$, by definition of $\pi_{[g]}$,
there is a $\gamma_g\in\Omega(G,e,g)$ such that
$[\gamma_g]=[\sigma\circ\gamma_g]$ in $\pi_1(G,e,g)$. That is,
$[\gamma_g]\in(\widetilde{G})^{\widetilde{\sigma}}$. By Lemma
\ref{L:finite}, $(\widetilde{G})^{\widetilde{\sigma}}$ is
connected. So
$g=\pi([\gamma_g])\in\pi((\widetilde{G})^{\widetilde{\sigma}})=G^\sigma_0$,
that is, $[g]=G^\sigma_0$.
\end{proof}

\begin{proposition}\label{P}
Each nontrivial element of $\Pi^\sigma$ has order $2$, and the
homomorphism $f:G^\sigma/G^\sigma_0\rightarrow\Pi^\sigma$ is an
isomorphism.
\end{proposition}

\begin{proof}
The first assertion follows from Lemma \ref{L:inverseclosed}, the
surjectivity of $f$ follows from the definition of $\Pi^\sigma$,
and the injectivity of $f$ follows from Lemma \ref{L:injective}.
\end{proof}

In particular, we have

\begin{lemma}\label{L:disjoint}
For $g_1,g_2\in G^\sigma$, if $[g_1]\neq[g_2]$, then
$\pi_{[g_1]}\cap\pi_{[g_2]}=\emptyset$. \qed
\end{lemma}

Now we are prepared to prove Theorem \ref{T:covering} and Theorem
\ref{T:quotient}.

\vskip 0.3cm {\flushleft {\it Proof of Theorem \ref{T:quotient}.}}
By Lemma \ref{L:finite}, $G^\sigma/G^\sigma_0$ is a finite group.
By Proposition \ref{P}, $G^\sigma/G^\sigma_0$ is abelian, and each
nontrivial element of $G^\sigma/G^\sigma_0$ has order $2$. By the
structure theorem of finitely generated abelian groups,
$G^\sigma/G^\sigma_0$ is isomorphic to $(\mathbb{Z}_2)^r$ for some
non-negative integer $r$. If $\pi_1(G)$ is finite with odd order,
so are $\pi_{(G^\sigma)}$, $\Pi^\sigma$, and
$G^\sigma/G^\sigma_0$. This forces $r=0$.\qed

\vskip 0.3cm {\flushleft {\it Proof of Theorem \ref{T:covering}.}}
We define an equivalence relation on the set $\bigcup_{g\in
G}\Omega(G,e,g)$ as follows. For
$\gamma_1,\gamma_2\in\bigcup_{g\in G}\Omega(G,e,g)$,
$\gamma_1\sim\gamma_2$ if and only if $\gamma_1(1)=\gamma_2(1)$
and $[\gamma_{12}]\in\pi_K$, where $\gamma_{12}\in\Omega(G,e,e)$
is defined by
$$
\gamma_{12}(t)=
\begin{cases}
\gamma_1(2t) & t\in[0,\frac{1}{2}];\\
\gamma_2(2-2t) & t\in[\frac{1}{2},1].
\end{cases}
$$
Denote $G'=(\bigcup_{g\in G}\Omega(G,e,g))/\sim$, and denote the
equivalence class of a $\gamma\in\bigcup_{g\in G}\Omega(G,e,g)$ by
$\langle\gamma\rangle$. Endow $G'$ with a smooth manifold
structure in the canonical way, and define the group operation on
$G'$ as follows. For
$\langle\gamma_1\rangle,\langle\gamma_2\rangle\in G'$, let
$\alpha(t)=\gamma_1(t)\gamma_2(t)$ and
$\beta(t)=\gamma_1(t)^{-1}$, and define
$\langle\gamma_1\rangle\langle\gamma_2\rangle=\langle\alpha\rangle,
\langle\gamma_1\rangle^{-1}=\langle\beta\rangle$. It is easy to
check that these are well defined, making $G'$ a Lie group. Let
$\pi:G'\rightarrow G$ be the map
$\pi(\langle\gamma\rangle)=\gamma(1)$, then $\pi$ is a covering
homomorphism. Since $\sigma_*(\pi_K)=\pi_K$, the map
$\sigma':G'\rightarrow G',
\sigma'(\langle\gamma\rangle)=\langle\sigma\circ\gamma\rangle$ is
a well defined involution of $G'$. It is obvious that
$\pi\circ\sigma'=\sigma\circ\pi$, that is, $\sigma'$ is just the
induced involution of $\sigma$ on $G'$. Now we have
\begin{align*}
&\langle\gamma\rangle\in(G')^{\sigma'}\\
\Leftrightarrow
&\langle\sigma\circ\gamma\rangle=\langle\gamma\rangle\\
\Leftrightarrow &\sigma\circ\gamma\sim\gamma\\
\Leftrightarrow &\gamma(1)\in G^\sigma,
[\gamma_{\gamma(1)}^\sigma]\in\pi_K\\
\Leftrightarrow &\gamma(1)\in K \quad (\text{by Lemma
\ref{L:disjoint}})\\
\Leftrightarrow &\langle\gamma\rangle\in\pi^{-1}(K).
\end{align*}
That is, $(G')^{\sigma'}=\pi^{-1}(K)$. This completes the proof of
Theorem \ref{T:covering}.\qed

\begin{remark}
The covering group $G'$ of $G$ satisfying the conclusion of
Theorem \ref{T:covering} is not necessarily unique. For example,
let $G$ be a semisimple Lie group with trivial center, and let
$\sigma$ be a global Cartan involution of $G$. Then $K=G^\sigma$
is a maximal compact subgroup of $G$. Let $G'$ be any finite cover
of $G$ with an involution $\sigma'$ such that
$\pi\circ\sigma'=\sigma\circ\pi$, then $(G')^{\sigma'}$ is a
maximal compact subgroup of $G'$. Since $\pi^{-1}(G^\sigma)$ is a
compact subgroup of $G'$ containing $(G')^{\sigma'}$, it must
equal $(G')^{\sigma'}$. In fact, using the same method as in the
proof of Theorem \ref{T:covering}, it is easy to see that a
covering group $G'$ of $G$ with covering homomorphism $\pi$
satisfying $\pi^{-1}(K)=(G')^{\sigma'}$ if and only if
$\pi_*(\pi_1(G',e))\cap\pi_{(G^\sigma)}=\pi_K$, and the group $G'$
that we constructed in the proof of Theorem \ref{T:covering} is
just the one satisfying $\pi_*(\pi_1(G',e))=\pi_K$.
\end{remark}

\begin{remark}
Professor Jiu-Kang Yu pointed out to the first author proofs
of Theorems 2.1 and 2.2 using nonabelian cohomology after he read
the first draft of the paper.
\end{remark}


\vskip 1.0cm
\section{Covering groups and the realization of symmetric spaces}

\vskip 0.3cm Let $G$ be a connected Lie group with an involution
$\sigma$, and let $K$ be an open subgroup of $G^\sigma$. Let
$\g=\ak\oplus\ap$ be the canonical decomposition of the Lie
algebra $\g$ of $G$ with respect to the differential of the
involution $\sigma$. Suppose $Ad_G(K)|_\ap$ is compact, that is,
$G/K$ has a structure of Riemannian symmetric space. Suppose
$\pi:G'\rightarrow G$ is a covering homomorphism of Lie groups.
Then the twisted conjugate action $\tau$ of $G$ on $G$, which is
defined by $\tau_g(h)=gh\sigma(g)^{-1}$, can be lifted to a well
defined action $\tau'$ of $G$ on $G'$, that is,
$\tau'_g(h')=g'h'\sigma'(g')^{-1}, g\in G, h'\in G'$, where
$g'\in\pi^{-1}(g)$, $\sigma'$ is the induced involution of
$\sigma$ on $G'$.

\begin{theorem}\label{T:realization}
Under the above assumptions, there is a covering group $G'$ of $G$
with covering homomorphism $\pi$ such that $P'=\exp_G'(\ap)$ is a
closed submanifold of $G'$, and such that the map
$\varphi:G/K\rightarrow P'$ defined by
$\varphi(gK)=g'\sigma'(g')^{-1}$ is a diffeomorphism, where
$g'\in\pi^{-1}(g)$, $\sigma'$ is the induced involution of
$\sigma$ on $G'$. With respect to the actions of $G$ by left
multiplication on $G/K$ and by the action $\tau'$ on $P'$,
$\varphi$ is equivariant.
\end{theorem}

\begin{proof}
By Theorem \ref{T:covering}, there is a covering group $G'$ of $G$
with covering homomorphism $\pi$ such that
$\pi^{-1}(K)=(G')^{\sigma'}$. So $\pi$ induces a $G$-equivariant
diffeomorphism $\overline{\pi}:G'/(G')^{\sigma'}\rightarrow G/K$
with respect to left multiplications (note that the left
multiplication of $G$ on $G'/(G')^{\sigma'}$ is well defined).
Since $Ad_{G'}((G')^{\sigma'})|_\ap=Ad_G(K)|_\ap$ is compact, by
Corollary 2.6 in \cite{AW}, $P'=\exp_{G'}(\ap)$ is a closed
submanifold of $G'$, and the map
$\varphi':G'/(G')^{\sigma'}\rightarrow P',
\varphi'(g'(G')^{\sigma'})=g'\sigma'(g')^{-1}$ is a
$G'$-equivariant diffeomorphism with respect to left
multiplication and twisted conjugate action of $G'$. Let
$\varphi=\varphi'\circ(\overline{\pi})^{-1}$. Then
$\varphi(gK)=\varphi'(g'(G')^{\sigma'})=g'\sigma'(g')^{-1}$. It is
obviously a $G$-equivariant diffeomorphism with respect to left
multiplication on $G/K$ and the action $\tau'$ on $P'$.
\end{proof}

\end{document}